\documentclass[12pt]{amsart} 
\usepackage{amsfonts,graphics,amsmath,amsthm,amsfonts,amscd, amssymb,amsmath,latexsym,multicol}
\usepackage{epsfig}
\usepackage{flafter}

\input diagrams

\makeatletter

\def\jobis#1{FF\fi
  \def\predicate{#1}%
  \edef\predicate{\expandafter\strip@prefix\meaning\predicate}%
  \edef\job{\jobname}%
  \ifx\job\predicate
}

\makeatother

\if\jobis{proposal}%
 \def\try{subsection}%
\else
  \def\try{section}%
\fi

\theoremstyle{plain}
\newtheorem{theorem}{Theorem}[\try]
\newtheorem{corollary}[theorem]{Corollary}
\newtheorem{lemma}[theorem]{Lemma}

\newtheorem{proposition}[theorem]{Proposition}
\newtheorem{definition-lemma}[theorem]{Definition-Lemma}
\newtheorem{question}[theorem]{Question}
\newtheorem{definition}[theorem]{Definition}
\newtheorem{remark}[theorem]{Remark}

\def\ideal#1.{I_{#1}}
\def\ring#1.{\mathcal {O}_{#1}}
\def\fring#1.{\hat{\mathcal {O}}_{#1}}
\def\proj#1.{\mathbb {P}(#1)}
\def\pr #1.{\mathbb {P}^{#1}}
\def\dpr #1.{\hat{\mathbb {P}}^{#1}}
\def\af #1.{\mathbb A^{#1}}
\def\Hz #1.{\mathbb F_{#1}}
\def\Hbz #1.{\overline{\mathbb F}_{#1}}
\def\fb#1.{\underset #1 {\times}}
\def\rest#1.{\underset {\ \ring #1.} \to \otimes}
\def\au#1.{\operatorname {Aut}\,(#1)}
\def\deg#1.{\operatorname {deg } (#1)}
\def\pic#1.{\operatorname {Pic}\,(#1)}
\def\pico#1.{\operatorname{Pic}^0(#1)}
\def\picg#1.{\operatorname {Pic}^G(#1)}
\def\ner#1.{NS (#1)}
\def\rdown#1.{\llcorner#1\lrcorner}
\def\rfdown#1.{\lfloor{#1}\rfloor}
\def\rup#1.{\ulcorner{#1}\urcorner}
\def\rcup#1.{\lceil{#1}\rceil}
\def\cone#1.{\operatorname {NE}(#1)}
\def\ccone#1.{\overline{\operatorname {NE}}(#1)}
\def\coef#1.{\frac{(#1-1)}{#1}}
\def\vit#1.{D_{\langle #1 \rangle}}
\def\mm#1.{\overline {M}_{0,#1}}
\def\H1#1.{H^1(#1,{\ring #1.})}
\def\ac#1.{\overline {\mathbb F}_{#1}}

\def\adj#1.{\frac {#1-1}{#1}}
\def\spn#1.{\overline{#1}}
\def\pek#1.#2.{\Cal P^{#1}(#2)}
\def\plk#1.#2.{\Cal P^{\leq #1}(#2)}
\def\ev#1.{\operatorname{ev_{#1}}}
\def\ilist#1.{{#1}_1,{#1}_2,\dots}
\def\bminv#1.{(\nu_1,s_1;\nu_2,s_2;\dots ;\nu_{#1},s_{#1};\nu_{r+1})}
\def\zinv#1.{(\nu_1,s_1;\nu_2,s_2;\dots ;\nu_{#1},s_{#1};0)}
\def\iinv#1.{(\nu_1,s_1;\nu_2,s_2;\dots ;\nu_{#1},s_{#1};\infty)}

\def\llist#1.#2.{{#1}_1,{#1}_2,\dots,{#1}_{#2}}
\def\lomitlist#1.#2.{{#1}_1,{#1}_2,\dots,\hat {{#1}_i}, \dots, {#1}_{#2}}
\def\lomitlistz#1.#2.{{#1}_0,{#1}_1,\dots,\hat {{#1}_i}, \dots, {#1}_{#2}}
\def\loc#1.#2.{\Cal O_{#1,#2}}
\def\fderiv#1.#2.{\frac {\partial #1}{\partial #2}}
\def\deriv#1.#2.{\frac {d #1}{d #2}}
\def\map#1.#2.{#1 \longrightarrow #2}
\def\rmap#1.#2.{#1 \dasharrow #2}
\def\emb#1.#2.{#1 \hookrightarrow #2}
\def\non#1.#2.{\text {Spec }#1[\epsilon]/(\epsilon)^{#2}}
\def\Hi#1.#2.{\text {Hilb}^{#1}(#2)}
\def\sym#1.#2.{\operatorname {Sym}^{#1}(#2)}
\def\Hb#1.#2.{\text {Hilb}_{#1}(#2)}
\def\Hm#1.#2.{\Hom_{#1}(#2)}
\def\prd#1.#2.{{#1}_1\cdot {#1}_2\cdots {#1}_{#2}}
\def\Bl #1.#2.{\operatorname {Bl}_{#1}#2}
\def\pl #1.#2.{#1^{\otimes #2}}
\def\mgn#1.#2.{\overline {M}_{#1,#2}}
\def\ialist#1.#2.{{#1}_1 #2 {#1}_2, #2\dots}
\def\pair#1.#2.{\langle #1, #2\rangle}
\def\vandermonde#1.#2.{\left|
\begin{matrix}
1 & 1 & 1 & \dots & 1\\
{#1}_1 & {#1}_2 & {#1}_3 & \dots & {#1}_{#2}\\
{#1}_1^2 & {#1}_2^2 & {#1}_3^2 & \dots & {#1}_{#2}^2\\
\vdots & \vdots & \vdots & \ddots & \vdots\\
{#1}_1^{#2-1} & {#1}_2^{#2-1} & {#1}_2^{#2-1} & \dots & {#1}_{#2}^{#2-1}\\
\end{matrix}
\right|
}
\def\vandermondet#1.#2.{\left|
\begin{matrix}
1 & {#1}_1   & {#1}_1^2 & \dots & {#1}_1^{#2-1}\\
1 & {#1}_2   & {#1}_2^2 & \dots & {#1}_2^{#2-1}\\
1 & {#1}_3   & {#1}_3^2 & \dots & {#1}_3^{#2-1}\\
\vdots & \vdots & \vdots & \ddots & \vdots\\
1 & {#1}_{#2}& {#1}_{#2}^2 & \dots & {#1}_{#2}^{#2-1}\\
\end{matrix}
\right|
}
\def\gr#1.#2.{\mathbb{G}(#1,#2)}

\def\alist#1.#2.#3.{{#1}_1 #2 {#1}_2 #2\dots #2 {#1}_{#3}}
\def\zlist#1.#2.#3.{#1_0 #2 #1_1 #2\dots #2 #1_{#3}}
\def\lomitlist30#1.#2.#3.{{#1}_0,{#1}_1 #2 \dots #2\hat {{#1}_i} #2\dots #2 {#1}_{#3}}
\def\lmap#1.#2.#3.{#1 \overset{#2}{\longrightarrow} #3}
\def\mes#1.#2.#3.{#1 \longrightarrow #2 \longrightarrow #3}
\def\ses#1.#2.#3.{0\longrightarrow #1 \longrightarrow #2 \longrightarrow #3 \longrightarrow 0}
\def\les#1.#2.#3.{0\longrightarrow #1 \longrightarrow #2 \longrightarrow #3}
\def\res#1.#2.#3.{#1 \longrightarrow #2 \longrightarrow #3\longrightarrow 0}
\def\Hi#1.#2.#3.{\text {Hilb}^{#1}_{#2}(#3)}
\def\ten#1.#2.#3.{#1\underset {#2}{\otimes} #3}
\def\lomitlist30#1.#2.#3.{{#1}_0 #2 {#1}_1 #2 \dots #2 \hat {{#1}_i} #2 \dots #2 {#1}_{#3}}

\def\Hom{\operatorname{Hom}}

\def\dim{\operatorname{dim}}

\def\deg{\operatorname{deg}}

\def\lcs{\operatorname{LCS}}

\def\rest{\operatorname{res}}

\def\ch{\operatorname{CH}}

\def\e{\Cal E}

\def\e1{E_1}
\def\e2{E_2}

\def\mapdown#1{\big\downarrow\rlap{$\vcenter
{\hbox{$\scriptstyle#1$}}$}}

\def\mapse#1{
{\vcenter{\hbox{$\mathop{\smash{\raise1pt\hbox{$\diagdown$}\!\lower7pt
\hbox{$\searrow$}}\vphantom{p}}\limits_{#1}\vphantom{\mapdown{}}$}}}}

\def\VR#1.{height#1pt&\omit&&\omit&&\omit&&\omit&&\omit&\cr}

\def\VRT#1.{height#1pt&\omit&&\omit&\cr}

\begin{document}
\title{Shokurov's Rational Connectedness Conjecture}
\author{Christopher D. Hacon} 
\address{Department of Mathematics \\  
University of Utah\\  
155 South 1400 E\\
JWB 233\\
Salt Lake City, UT 84112, USA}
\email{hacon@math.utah.edu}
\author{James M\textsuperscript{c}Kernan} 
\address{Department of Mathematics\\ 
University of California at Santa Barbara\\ 
Santa Barbara, CA 93106, USA} 
\email{mckernan@math.ucsb.edu}

\begin{abstract} We prove a conjecture of V. V. Shokurov which in particular implies 
that the fibers of a resolution of a variety with divisorial log terminal singularities
are rationally chain connected.
\end{abstract} 

\thanks{The first author was partially supported by NSA research grant no:
  MDA904-03-1-0101 and by a grant from the Sloan Foundation.  Part of this work was
  completed whilst the second author was visiting the Pontificia Universidad Cat\'olica
  del Per\'u, Lima and the Tokyo Mathematics Department.}

\maketitle
\pagestyle{plain}
\section{Introduction}
\label{s_introduction} 

In recent years it has become increasingly clear that the geometry of higher dimensional
varieties is closely related to the geometry of rational curves on these varieties. From
the point of view of the minimal model program, one expects that rational curves on
varieties with mild singularities (e.g. log terminal singularities) share many of the
basic properties of rational curves on smooth varieties.  Surprisingly, very little seems
to be known in this direction.  The purpose of this paper is to give an affirmative answer
to several natural questions that arise in this context.  For example we show that:

\textit{If $(X, \Delta)$ is a divisorially log terminal pair and
  $f:\map Y.X.$ is a birational morphism, then for any $x\in X$,
  $f^{-1}(x)$ is rationally chain connected and in particular covered
  by rational curves.}

Two immediate consequences of this are:
\begin{enumerate} 
\item \textit{If $(X, \Delta )$ is a divisorially log terminal pair, and $g: \rmap X.Z.$
  is a rational map to a proper variety which is not everywhere defined, then $Z$ contains
  a rational curve.}
\item \textit{If $(X, \Delta )$ is a divisorially log terminal pair,
  then $X$ is rationally chain connected if and only if it is
  rationally connected.}
\end{enumerate}  

We now turn to a more detailed discussion of the results of this paper.

\begin{definition}\label{d_mod} Let $X$ be a reduced, separated scheme of finite type 
over an algebraically closed field (so that every irreducible component of $X$ is a
variety) and let $V$ be any subset.  We will say that a curve $C$ is a \textbf{chain
  modulo $V$} if $C$ union a subset of $V$ is connected.  We will say that \textbf{$X$ is
  rationally chain connected modulo $V$}, if any two points $x$ and $y$, which belong to
the same connected component of $X$, belong to a chain of rational curves modulo $V$.
\end{definition}

Note that if $V$ is empty and $X$ is irreducible, then $X$ is rationally chain connected
in the usual sense.  Note also that the disjoint union of two copies of $\pr 1.$ is
rationally chain connected but not connected.  Here is the main result of this paper:

\begin{theorem}\label{t_main} Let $(X,\Delta)$ be a log pair and let $f\colon\map X.S.$ be 
a projective morphism such that $-K_X$ is relatively big and $\ring X.(-m(K_X+\Delta))$ is
relatively generated, for some $m>0$.  Let $g\colon\map Y.X.$ be any birational morphism
and let $\pi\colon\map Y.S.$ be the composite morphism.

Then every fibre of $\pi$ is rationally chain connected modulo the inverse image of the
locus of log canonical singularities.
\end{theorem}

We note that we can weaken the hypothesis if we assume that $K_X+\Delta$ is kawamata log
terminal, see \eqref{l_semi}.  This result has some interesting consequences.  We start
with a very general result about the fundamental group of a log pair:

\begin{corollary}\label{c_main} Let $(X,\Delta)$ be a projective log pair such that 
$-(K_X+\Delta)$ is semiample and $-(K_X+\Delta)$ is big.

Then the fundamental group of $X$ is a quotient of the fundamental group of the locus of log
canonical singularities of the pair $(X,\Delta)$.  
\end{corollary}

Some other interesting consequences of \eqref{t_main} arise when we eliminate the locus of
log canonical singularities from the statement:

\begin{corollary}\label{c_rcck} Let $(X,\Delta)$ be a kawamata log terminal pair. 
Let $f\colon\map X.S.$ be a projective morphism such that $-(K_X+\Delta)$ is relatively nef 
and $-K_X$ is relatively big.  

 Then every fibre of $f$ is rationally chain connected. 
\end{corollary}

 With a little more work, we can also prove:

\begin{corollary}\label{c_rccl} Let $(X,\Delta)$ be a log canonical pair.  Let 
$f\colon\map X.S.$ be a projective morphism such that $-(K_X+\Delta)$ is relatively ample.

 Then every fibre of $f$ is rationally chain connected. 
\end{corollary}

\eqref{c_rccl} was conjectured by Shokurov in \cite{Shokurov00a}, and was proved there,
assuming the MMP.  In the same paper, he also conjectured the following, which he proved under
the same assumptions:

\begin{corollary}\label{c_dlt} Let $(X,\Delta)$ be a divisorially log terminal pair. 

If $g\colon\map Y.X.$ is any birational morphism then the fibres of $g$ are rationally chain
connected.
\end{corollary}

It was pointed out in \cite{Shokurov00a} that one then gets the following:

\begin{corollary}\label{c_covered} Let $f:\rmap X.Y.$ be a rational morphism of normal proper 
varieties such that $(X,\Delta)$ is a divisorially log terminal pair for some effective
divisor $\Delta$.  Then, for each closed point $x\in X$, the indeterminacy locus of $x$ is
covered by rational curves.
\end{corollary}

Recall that if $W\subset X\times Y$ is the closure of the graph of $f$ and $p$ and $q$ are
the projections from $W$ to $X$, and $Y$, then the indeterminacy locus of $x$ is defined
as $q(p^{-1}(x))$.

Another very interesting consequence of \eqref{c_dlt} is:

\begin{corollary}\label{c_equivalent} Let $(X,\Delta)$ be a divisorially log terminal pair.  
Then $X$ is rationally chain connected iff it is rationally connected.  
\end{corollary}

 The same methods which are used to prove \eqref{t_main}, in conjunction with 
\eqref{c_equivalent}, yield the following:

\begin{corollary}\label{c_over} Let $(X,\Delta)$ be a kawamata log terminal pair.  
Let $f\colon\map X.S.$ be a projective morphism with connected fibres such that
$-(K_X+\Delta)$ is relatively nef and $-K_X$ is relatively big.  Let $g\colon\map Y.X.$ be
any birational morphism, and let $\pi\colon\map Y.S.$ be the composition.  Let $T$ be any
irreducible subset of $S$ and let $W$ be the inverse image of $T$ inside $Y$.

Then there is an irreducible closed subset $E$ of $W$, which dominates $T$ and which has
connected and rationally connected fibres.
\end{corollary}

One interesting feature of rationally connected varieties is that a family of rationally
connected varieties over a curve always admits a section.  As a consequence of
\eqref{c_over}, we are able to show that a similar result holds for families of Fano
varieties, a result which was also conjectured by Shokurov:

\begin{corollary}\label{c_section} Let $(X,\Delta)$ be a kawamata log terminal pair.  
Let $f\colon\map X.S.$ be a projective morphism with connected fibres such that
$-(K_X+\Delta)$ is relatively nef and $-K_X$ is relatively big.  Let $g\colon\map Y.X.$ be
any birational morphism, and let $\pi\colon\map Y.S.$ be the composition.

 Then $\pi$ has a section over any curve.  
\end{corollary}

Another reason why rationally connected and rationally chain connected varieties are
interesting, is because their topology is particularly simple.  In particular this means
that the intersection theory of a divisorially log terminal pair and a resolution are closely
related:

\begin{corollary}\label{c_chow} Let $(X,\Delta)$ be a kawamata log terminal pair.  
Let $f\colon\map X.S.$ be a projective morphism with connected fibres such that
$-(K_X+\Delta)$ is relatively nef and $-K_X$ is relatively big.  Let $g\colon\map Y.X.$ be
any birational morphism, and let $\pi\colon\map Y.S.$ be the composition.

Then the natural map 
$$
\pi_*\colon\map \ch^0(Y).\ch^0(S).,
$$
is an isomorphism.  
\end{corollary}

   The following, although not a direct consequence of \eqref{t_main}, is proved using very
similar methods:

\begin{corollary}\label{c_iitaka} Let $(X,\Delta)$ be a kawamata log terminal projective pair
and let $\rmap X.Z.$ be the Iitaka fibration associated to $-(K_X+\Delta)$.

 Then $Z$ is rationally connected.  
\end{corollary}

 As a consequence we obtain:
\begin{corollary}\label{c_fano} Let $(X,\Delta)$ be a kawamata log terminal projective pair
and suppose that that $-(K_X+\Delta)$ is big and nef.  

 Then $X$ is rationally connected.  
\end{corollary}

We remark that \eqref{c_fano} has also been recently proved by Zhang, using a similar
argument, see \cite{Zhang04}.  We were working on this paper, when his result appeared on
the archive.  As Zhang points out in \cite{Zhang04}, \eqref{c_fano} implies the following:

\begin{corollary}\label{c_simply} Let $(X,\Delta)$ be a kawamata log terminal projective pair
and suppose that that $-(K_X+\Delta)$ is big and nef. 

 Then $X$ is simply connected.  
\end{corollary}

It should be pointed out though that it is more natural to prove the stronger result that
the smooth locus of $X$ has finite fundamental group.  The only known case of this much
stronger result, is for surfaces, see \cite{KM99}.  It is proved in \cite{McKernan02} that
at least the algebraic fundamental group of the smooth locus is finite.

\begin{remark}\label{r_replace} Note that \eqref{c_rcck} and \eqref{c_simply}, which 
are stated for kawamata log terminal pairs extend to the case when the locus of log
canonical singularities is a point (equivalently, by connectedness, when $V$ has dimension
zero).
\end{remark}

We now give a quick sketch of the proof of \eqref{t_main}.  Let us suppose that we want to
prove that $F$ is rationally connected.  By the main result of \cite{GHS03}, it suffices
to prove that whenever we have a test rational map $t\colon\map F.Z.$, then either $Z$ is
uniruled or it is a point, see \eqref{l_sufficient}.  Suppose $Z$ is not uniruled.  By the
main result of \cite{BDPP04}, $K_Z$ is pseudo-effective.

Suppose for a moment, that we can produce a divisor $\Theta$ with three key properties:
\begin{itemize} 
\item[(a)] $K_F+\Theta$ has Kodaira dimension zero, 
\item[(b)] there is an ample $\mathbb{Q}$-divisor $H$ on $Z$ such that $t^*H\leq \Theta$, and 
\item[(c)] $K_F+\Theta$ is log terminal on the general fibre.  
\end{itemize} 

By log additivity, it follows that the Kodaira dimension of $K_F+\Theta$ is at least the
dimension of $Z$, so that $Z$ is a point by (a), see \eqref{p_uni}.  It remains to
indicate how to produce the divisor $\Theta$.  In our case, $F$ is a smooth divisor and a
component of the fibre of a resolution of a Fano fibration $f\colon\map X.S.$.  The
existence of a divisor $\Delta$ on the total space of the Fano fibration satisfying (a),
(b) and (c) is quite straightforward, and it is equally straightforward to lift this to a
resolution $g\colon\map Y.X.$ of the total space.  Conditions (b) and (c) then descend to
$F$.  The tricky part is ensuring that condition (a) continues to hold.

The first step is to realise $F$ as a component of the locus of log canonical
singularities, with respect to some divisor.  At this point we apply the main technical
result of \cite{HM05b}, which says, roughly speaking, that we can lift appropriate sections from $F$
to $Y$.

\section{Notation and conventions}
\label{s_notation}

We work over the field of complex numbers $\mathbb{C}$.  A $\mathbb{Q}$-Cartier divisor
$D$ on a normal variety $X$ is \textit{nef} if $D\cdot C\geq 0$ for any curve $C\subset
X$.  We say that two $\mathbb{Q}$-divisors $D_1$, $D_2$ are $\mathbb{Q}$-linearly
equivalent ($D_1\sim _{\mathbb{Q}} D_2$) if there exists an integer $m>0$ such that $mD_i$
are linearly equivalent.  Given a morphism of normal varieties $g\colon\map Y.X.$, we say
that two $\mathbb{Q}$-divisors $D_1$ and $D_2$ are $\mathbb{Q}$-$g$-linearly equivalent
($D_1\sim _{g,\mathbb{Q}} D_2$ if there is a positive integer $m$ and a Cartier divisor
$B$ on $X$ such that $mD_1\sim mD_2+f^*B$.  We say that a $\mathbb{Q}$-Weil divisor $D$ is
big if we may find an ample divisor $A$ and an effective divisor $B$, such that $D \sim
_{\mathbb{Q}} A+B$.  A \textit{log pair} $(X,\Delta)$ is a normal variety $X$ and an
effective $\mathbb{Q}$-Weil divisor $\Delta$ such that $K_X+\Delta$ is
$\mathbb{Q}$-Cartier.  A projective morphism $g \colon\map Y.X.$ is a \textit{log
  resolution} of the pair $(X,\Delta )$ if $Y$ is smooth and $g^{-1}(\Delta )\cup
\{\,\text {exceptional set of $g$}\,\}$ is a divisor with normal crossings support.  We
write $g^*(K_X +\Delta )=K_Y +\Gamma$ and $\Gamma =\sum a_i\Gamma _i$ where $\Gamma _i$
are distinct reduced irreducible divisors.  The log discrepancy of $\Gamma_i$ is $1-a_i$.
The \textit{locus of log canonical singularities of the pair $(X,\Delta)$}, denoted
$\lcs(X,\Delta)$, is equal to the image of those components of $\Gamma$ of coefficient at
least one (equivalently log discrepancy at most zero).  The pair $(X,\Delta )$ is
\textit{kawamata log terminal} if for every (equivalently for one) log resolution
$g\colon\map Y.X.$ as above, the coefficients of $\Gamma$ are strictly less than one, that
is $a_i<1$ for all $i$.  Equivalently, the pair $(X,\Delta)$ is kawamata log terminal if
the locus of log canonical singularities is empty.  We say that the pair $(X,\Delta)$ is
\textit{divisorially log terminal} if the coefficients of $\Delta$ lie between zero and
one, and there is a log resolution such that the coefficients of the $g$-exceptional
divisors are all less than one.

Given a pair $(X,\Delta)$ and a birational morphism, we will often consider decompositions
of the form
$$
K_Y+\Gamma=\pi^*(K_X+\Delta)+E,
$$
where $\Gamma$ and $E$ are effective, often with no common components, and where $E$ is
exceptional.  Here, even though the equals sign denotes $\mathbb{Q}$-linearly equivalence,
in fact it is often implicitly assumed that $\Gamma$ is a combination of exceptional
divisors and the strict transform of $\Delta$.  In fact, under the further assumption that
$\Gamma$ and $E$ have no common components, such a decomposition is unique.  However, we
will often make no such assumption concerning the support of $\Gamma$, so that such a
decomposition is far from unique.

Let $\phi\colon\rmap X.Y.$ be a rational map.  We say that a subset $V$ of $X$
\textit{dominates} $Y$, if the inverse image of $V$, in the graph of $\phi$, dominates
$Y$.

Let $f\colon\map X.S.$ be a morphism.  We say that a divisor $E$ is exceptional for $f$,
if $f(E)$ has codimension at least two in $S$.

\section{Some Examples}
\label{s_examples}

In this section we collect together some examples, whose purpose is to show that the
results stated in \S \ref{s_introduction} are in some sense best possible, and to motivate
their proofs.

We remark that in theorem \eqref{c_fano}, one can not remove the hypothesis that
$-(K_X+\Delta)$ is nef or that $(X,\Delta )$ is kawamata log terminal as shown by the
following well known example.  Let $f\colon\map S.C.$ be any $\pr 1.$-bundle over an
elliptic curve and let $E$ be a section of minimal self-intersection.  Then $f$ is the
maximal rationally connected fibration of $S$.  Suppose first that $E^2<0$.  Then
$$
-(K_S+tE) \qquad \text{is} \qquad  
\begin{cases} \text{big}   & \text{for any $t<2$} \\
              \text{ample} & \text{for any $1<t<2$} \\
              \text{nef}   & \text{for any $1\leq t\leq 2$} \\
              \text{lc}    & \text{for any $t\leq 1$} \\
              \text{klt}   & \text{for any $t<1$} \\
\end{cases}
$$
Clearly $S$ is not rationally chain connected.  If $t=1$ then $K_S+E$ is log canonical
and $-(K_S+E)$ is nef and big.  If $t<1$ then $K_S+tE$ is kawamata log terminal and
$-(K_S+tE)$ is big but not nef.  If $1<t<2$ then $-(K_S+tE)$ is ample but not log canonical.
Notice that by contracting the negative section $E$, we get a rationally chain connected
surface $T$.  Therefore, rational chain connectedness is not a birational property of $S$;
of course if $f\colon\map X'.X.$ is a birational morphism of varieties and $X'$ is
rationally chain connected, then $X$ is also rationally chain connected.

  One might well ask the following:
\begin{question}\label{q_ask} Let $(X,\Delta)$ be a klt pair, $\Delta$ effective and
$-(K_X+\Delta)$ nef.  If $X\dashrightarrow W$ is the MRCC fibration then does it follow
that 
$$
\dim X\geq \dim W +\kappa (-(K_X+\Delta))\ ?
$$ 
\end{question}

Now suppose that we let $E_2$ be the unique non-split extension of $\ring C.$ by $\ring
C.$ and set $S$ to be the projectivisation of $E_2$.  Then $E$ has self-intersection zero.
$-K_S=2E$, so $-K_S$ is nef of numerical dimension one.  But in fact no multiple of $E$
moves.  Indeed if it did, then there would be an \'etale cover $\pi\colon\map C'.C.$
such that $\pi^*E_2$ splits.  It is easy to see that in characteristic zero this never
happens.  Thus the Kodaira dimension of $-K_S$ is zero and the sum of the Kodaira
dimension plus the dimension of $C$ is less than the dimension of $S$ so that the answer
to \eqref{q_ask} is no, in general.  

Finally, we give an easy example, to illustrate the fact that one cannot expect every
component of the exceptional locus of a resolution of a divisorially log terminal
singularity to be rationally connected.  In fact it is rarely the case that each
individual component is rationally chain connected.

Let $X$ be a smooth threefold.  Pick a point $x\in X$ and let $\map Y_1.X.$ blow up this
point.  Now let $\map Y.Y_1.$ blow up a cubic curve in the exceptional divisor, which is a
copy of $\pr 2.$.  Consider the birational morphism $\pi\colon\map Y.X.$.  Denote the
first exceptional divisor, a copy of $\pr 2.$, by $E$ and the other by $F$.  It is easy to
check that $-(K_Y+4/5E+1/5F)$ is relatively ample.  Note that $E\cup F$ is rationally
chain connected, but that $F$ is not even rationally chain connected.  Indeed $F$ is once
again a $\pr 1.$-bundle over an elliptic curve.

\section{Uniruled, rationally connected and rationally chain connected varieties}
\label{s_uni}

In this section we give sufficient conditions to ensure that a variety is either uniruled,
or rationally connected, or rationally chain connected.  First a criterion 
to ensure that a variety is uniruled:

\begin{proposition}\label{p_uni} Let $(X,\Delta)$ be a projective log pair and let $h\colon\map X.F.$ 
and $t\colon\rmap F.Z.$ be a morphism and a rational map, where $F$ and $Z$ are
projective, with the following properties:
\begin{enumerate} 
\item the locus of log canonical singularities of  $K_X+\Delta$ does not dominate $Z$, where
$\rmap X.Z.$ is the composition of $t$ and $h$,
\item $K_X+\Delta$ has Kodaira dimension at least zero on the general fibre of $\rmap
X.Z.$, (that is, if $g\colon\map Y.X.$ resolves the indeterminacy of $\rmap X.Z.$ then
$g^*(K_X+\Delta)$ has Kodaira dimension at least zero on the general fibre of the induced
morphism $\map Y.Z.$),
\item $K_X+\Delta$ has Kodaira dimension at most zero, and
\item there is an ample divisor $A$ on $F$ such that $h^*A\leq\Delta$.  
\end{enumerate} 
 
 Then either $Z$ is a point or it is uniruled.  
\end{proposition}
\begin{proof} Suppose that $Z$ is not uniruled.  Blowing up $Z$, we may assume that $Z$ is
smooth.  Let $g\colon\map Y.X.$ be a log resolution of $X$, such that the induced rational
map $\rmap Y.Z.$ is in fact a morphism $\psi\colon\map Y.Z.$ .  We may write
$$
K_Y+\Theta=g^*(K_X+\Delta)+E,
$$
where $\Theta$ and $E$ are effective, with no common components and $E$ is exceptional.  Now
set $\Gamma=\Theta+\epsilon E'$, where $\epsilon$ is a sufficiently small rational number,
and $E'$ is the support of the exceptional locus.  With this choice of $\Gamma$, the
Kodaira dimension of $K_Y+\Gamma$ is at least zero on the general fibre of $\psi$, and the
support of $\Gamma$ contains the full exceptional locus.

As $Z$ is not uniruled, it follows by the main result of \cite{BDPP04} that $K_Z$ is
pseudo-effective.  By (4), and our choice of $\Gamma$, $\Gamma$ contains the pullback of
an ample divisor from $F$.  Possibly replacing $\Gamma$ by a linearly equivalent divisor,
we may find an ample divisor $G$ on $Z$ such that $\psi^*G\leq\Gamma$.  Since $K_Y+\Gamma$
is log terminal on the general fibre of $\psi$, it follows by log additivity of the
Kodaira dimension, see Corollary 2.11 of \cite{HM05b}, that the Kodaira dimension of
$K_Y+\Gamma$ is at least the dimension of $Z$.  As the Kodaira dimension of $K_Y+\Gamma$
is at most the Kodaira dimension of $K_X+\Delta$, it follows that $Z$ is a point.
\end{proof}

 It is easy to use \eqref{p_uni} to prove that a variety is rationally chain connected:

\begin{lemma}\label{l_sufficient} Let $F$ be a normal variety.  
\begin{enumerate} 
\item $F$ is rationally connected iff for every non-constant dominant rational map $t\colon\rmap
F.Z.$, $Z$ is uniruled.
\item $F$ is rationally chain connected modulo $V$ iff for every non-constant dominant
rational map $t\colon\rmap F.Z.$, either $Z$ is uniruled or $V$ dominates $Z$.
\end{enumerate} 
\end{lemma}
\begin{proof} It is clear that the image of every rationally connected variety is rationally 
connected.  In particular the image of every rationally connected variety is either
uniruled or a point.  If $V$ does not dominate $Z$, then pick two general points $x$ and
$y$ of $F$.  Then they are connected by a chain of rational curves modulo $V$, and the
image of one of these curves must be a rational curve through the general point of $Z$.

To prove the reverse direction, first observe that $F$ is uniruled, by applying the basic
criterion to the identity map $\map F.F.$.  Let $F'$ be a smooth model of $F$, and let
$\rmap F'.Z.$ be the maximal rationally connected fibration of $F'$.  Blowing up, we may
assume that this map is a morphism.  Since $F$ is birational to $F'$, this induces a
rational map $t\colon\rmap F.Z.$.  By assumption either $Z$ is a point, or $Z$ is
uniruled, or $V$ dominates $Z$.

As $F$ is uniruled, $\map F'.Z.$ is a non-trivial fibration.  This morphism has the
defining property that if a rational curve $C$ meets a very general fibre $G$, then $C$ is
contained in $G$.  By the main result of \cite{GHS03}, we can lift any rational curve
which passes through the general point of $Z$ to $F$.  Thus $Z$ is not uniruled.  If $Z$
is a point, then $F$ is rationally connected.  The only other possibility is that $V$
dominates $Z$.  In this case $F$ is certainly rationally chain connected modulo $V$.
\end{proof}

\begin{corollary}\label{c_uni} Let $(X,\Delta)$ be a log pair and let $h\colon\map X.F.$
be a morphism.  Suppose that for every rational map $t\colon\rmap
F.Z.$, either the locus of log canonical singularities of $K_X+\Delta$
dominates $Z$, or (2-4) of \eqref{p_uni} hold.
 
Then $F$ is rationally chain connected modulo the image $R$ of the locus where $K_X+\Delta$
is not kawamata log terminal. 
\end{corollary}
\begin{proof} We are going to apply the criterion of \eqref{l_sufficient}.  Suppose we are
given a rational map $t\colon\rmap F.Z.$.  If $R$ dominates $Z$ there is nothing to prove.
Otherwise we just need to prove that $Z$ is either uniruled or a point and for this we
just need to observe that conditions (1-4) of \eqref{p_uni} hold.  \end{proof}

\section{A fine analysis of the fibres of a Fano fibration}
\label{s_fine}

In this section we state a detailed theorem about the fibres of a Fano fibration, which
despite being technical in nature, we expect will be of independent interest.  We fix some
notation which will hold throughout the section.  Let $(X,\Delta)$ be a log pair and let
$f\colon\map X.S.$ be a morphism such that $K_X+\Delta \sim _{\mathbb{Q},f} 0$ and
$\Delta$ is $f$-big.

Let $s\in S$ be any closed point and let $K$ be any effective $\mathbb{Q}$-Cartier divisor
on $S$.  Let $g\colon\map Y.X.$ be any birational morphism, such that the fibre of the
composite morphism $\pi\colon\map Y.S.$ over $s$ union the exceptional locus of $g$ union
the strict transform of $f^*K+\Delta$ is a divisor with normal crossings.  Let $\llist
F.k.$ be the components of the fibre $\pi^{-1}(s)$, which have log discrepancy greater
than zero, with respect to $K_X+\Delta$, and let $F$ be their union.  

Given $t\in [0,1]$, set $\Delta_t=\Delta+tf^*K$, $\Gamma_t=\Gamma+t\pi^*K$ (where $\Gamma$
is defined below) and let $V_t$ be the closure of
$$
\lcs(Y,\Gamma_t)-\lcs(Y,\Gamma).
$$
We let $\Gamma'_t$ denote the fractional part of  $\Gamma_t$.  

\begin{theorem}\label{t_fibres} With the notation above, we may pick $K$ and divisors 
$\Gamma$ and $E$ on $Y$, with the following properties:
\begin{enumerate} 
\item The equation
$$
K_Y+\Gamma=g^*(K_X+\Delta)+E,
$$
holds, where $\Gamma$ and $E$ are effective, with no common components and $E$ is
$g$-exceptional.  Moreover, we may write $\Gamma=A+B$, where $A$ is $\pi$-ample and $B$ is
effective.
\item $F=V_1$.  
\item Possibly relabelling and rescaling, we may assume that there are rational numbers,
$0=\zlist t.,.k.=1$, such that $V_i=\alist F.\cup.i.$, where, here and elsewhere, we adopt
the shorthand subscript $i$ in lieu of $t_i$.  For $\epsilon$ sufficiently small,
$V_{t_i-\epsilon}=V_{t_{i-1}}$.  Denote by $\Theta_i$ the restriction of $\Gamma'_i$ to
$F_i$.
\item Let $H$ be any $\pi$-ample divisor.  Then there is a constant $M$ such that 
$$
h^0(X,\ring X.(m(K_{F_i}+\Theta_i)+H))\leq M
$$
for all sufficiently divisible positive integers $m$, and any $1\leq i\leq k$. 
\item $F_i$ is rationally chain connected modulo $W_i=F_i\cap \lcs(Y,\Gamma_{i-1})$.  
\item If $K_Y+\Gamma$ is kawamata log terminal then $F_1$ is rationally connected.  
\item If $K_X+\Delta$ is kawamata log terminal then we may pick
$\Gamma$ so that $K_Y+\Gamma$ is kawamata log terminal.
\end{enumerate} 
\end{theorem}

We prove each part of \eqref{t_fibres} in a series of lemmas: 

\begin{lemma}\label{l_one} (1) of \eqref{t_fibres} holds.  
\end{lemma}
\begin{proof} We may write
$$
K_Y+\Gamma_1=g^*(K_X+\Delta)+E_1,
$$
where $\Gamma_1$ and $E_1$ are effective, with no common components and $E_1$ is
$g$-exceptional.  Now let $E'$ be the sum of all the exceptional divisors, taken with
coefficient one.  Set $\Gamma_2=\Gamma_1+\delta E'$, $E_2=E_1+\delta E'$, for some
positive rational number $\delta$.  If we choose $\delta$ small enough, then we do not
change the locus of log canonical singularities.  As $\Gamma_2$ and the strict transform
of $\Delta$ union the exceptional locus have the same support, it follows that $\Gamma_2$
is $\pi$-big.  It follows that we may write $\Gamma_2\sim_{\mathbb{Q},\pi} A+B$, where $A$ is
$\pi$-ample and $B$ is effective.  Let
$$
K_Y+\Gamma=g^*(K_X+\Delta)+E,
$$
be the decomposition obtained by cancelling like terms on both sides of
$$
K_Y+((1-\epsilon)\Gamma_2+\epsilon B)+\epsilon A=g^*(K_X+\Delta)+E_2.
$$
Thus $\Gamma$ and $E$ are effective with no common components, $E$ is exceptional and of
course $\epsilon A\leq\Gamma$ is relatively ample.
\end{proof}

\begin{lemma}\label{l_two} (2) of \eqref{t_fibres} holds.    
\end{lemma}
\begin{proof} If we pick $K$ sufficiently singular at $s$, then we may assume that
$F\subset V_1$.  Possibly replacing $K$ by a linearly equivalent $\mathbb{Q}$-divisor, we
may then assume that $F=V_1$.  \end{proof}

\begin{lemma}\label{l_three} (3) of \eqref{t_fibres} holds.    
\end{lemma}
\begin{proof} By (2), for every $0\leq t\leq 1$, $V_t$ is a union of components of $F$.  Let
$t_i$ be the smallest value of $t$, such that $F_i$ is a component of $\lcs(Y,\Gamma_t)$.
Possibly perturbing $K$, we may assume that $t_i\neq t_j$, if $i\neq j$ and so possibly
re-ordering and rescaling, we may assume that $0<\alist t.<.k.=1$. \end{proof}

\begin{lemma}\label{l_four} (4) of \eqref{t_fibres} holds.    
\end{lemma}
\begin{proof} This is the most technical part of \eqref{t_fibres}. 

We are certainly free to replace $H$ by a multiple.  Pick a divisor $G$ on $X$ such that
$g^*G\geq (\dim X+2)H$.  Then 
\begin{align*} 
g_*\ring Y.(m(K_Y+\Gamma_i)+(\dim X+2)H)  &\subset g_*\ring Y.(g^*G+mE)\\
                                          &=\ring X.(G),\\
\end{align*} 
for $m$ sufficiently divisible.  Therefore $\pi _*\ring Y.(m(K_Y+\Gamma_i)+(\dim
X+2)H)\subset f_*\ring X.(G)$.  On the other hand, by (3.17) of \cite{HM05b}, since
$\Theta_i$ does not contain $F_i$, we may lift any section of
$$
\pi_*(\ring X.(m(K_{F_i}+\Theta_i)+H))=H^0(X,\ring X.(m(K_{F_i}+\Theta_i)+H)),
$$
to $\pi_*(\ring Y.(m(K_Y+\Gamma_i)+(\dim X+2)H))$.
\end{proof}

\begin{lemma}\label{l_five} (5) of \eqref{t_fibres} holds.    
\end{lemma}
\begin{proof} We are going to apply \eqref{c_uni}.  Let $\psi\colon\rmap F_i.Z.$ be any
rational map, such that the locus where $K_{F_i}+(\Gamma_i-F_i)|_{F_i}$ is not kawamata
log terminal does not dominate $Z$.  We just need to check that conditions (2-4) of
\eqref{p_uni} hold.  Now
$$
K_Y+\Gamma_i=E_i,
$$
where $E_i$ is effective and exceptional, and $F_i$ is not contained in the support of
$E_i$.  It follows that the Kodaira dimension of $K_Y+\Gamma_i$ restricted to the general
fibre of $\psi$ is at least zero.  We may assume that $W_i$ does not dominate $Z$, so that
$\Gamma_i=\Gamma_i'$, on the general fibre of $\psi$.  It follows that $K_{F_i}+\Theta_i$
has Kodaira dimension at least zero on the general fibre of $\psi$.  Thus (2) holds.  (3)
is an immediate consequence of (4) of \eqref{t_fibres}.  (4) holds by (1) of
\eqref{t_fibres}.
\end{proof}

\begin{lemma}\label{l_six} (6) of \eqref{t_fibres} holds.    
\end{lemma}
\begin{proof} Since $W_0=\lcs(Y,\Gamma)$ is by assumption empty, it follows that 
$F_1$ is rationally chain connected.  As $F_1$ is smooth, it is in fact rationally
connected.
\end{proof}

\begin{lemma}\label{l_seven} (7) of \eqref{t_fibres} holds.    
\end{lemma}
\begin{proof} Clear from the construction of $\Gamma$ given in \eqref{l_one}.  
\end{proof}

\section{Proof of \eqref{t_main}}
\label{s_proof}

\begin{lemma}\label{l_alter} If $-(K_X+\Delta)$ is semiample and $-K_X$ is big, 
then we may find a divisor $\Delta'\geq\Delta$, such that some multiple of $K_X+\Delta'$
is linearly equivalent to zero, where $\Delta'$ contains an ample divisor, and where the locus
of log canonical singularities of $K_X+\Delta'$ is contained in the locus of log canonical
singularities of $K_X+\Delta$.  
\end{lemma}
\begin{proof} If we pick a general element $D\in |-m(K_X+\Delta)|$, where $m$ is sufficiently 
divisible, then $K_X+\Delta+D/m$ has the same locus of log canonical singularities as
$K_X+\Delta$.  Replacing $\Delta$ by $\Delta+D/m$, we may assume therefore that $\Delta$
is big and that some multiple of $K_X+\Delta$ is linearly equivalent to zero.  As $\Delta$
is big, we may write $\Delta\sim_{\mathbb{Q}}A+B$, where $A$ is ample and $B$ is
effective.  Let $\Delta'=(1-\epsilon)\Delta+\epsilon A+\epsilon B$.  Then for $\epsilon$
sufficiently small, the locus of log canonical singularities of $K_X+\Delta'$ is contained
in the locus of log canonical singularities of $K_X+\Delta$.
\end{proof}

\begin{proof}[Proof of \eqref{t_main}] Let $s\in S$.  This result is local over
$s\in S$.  Passing to an open neighbourhood of $s\in S$, we may as well assume that
$-(K_X+\Delta)$ is semiample.  By \eqref{l_alter} we may assume therefore that
$\Delta=A+B$, where $A$ is ample and $B$ is effective, and that some multiple of
$K_X+\Delta$ is linearly equivalent to zero.

Note that as $g$ has connected fibres, it follows that if $\map Y'.Y.$ is any birational
map, then we are free to replace $Y$ by $Y'$.  Thus we may assume that the hypothesis and
notation of \S \ref{s_fine} hold.

By induction on $i$, it suffices to prove that $F_i$ is rationally chain connected modulo
$W_i$, which is (5) of \eqref{t_fibres}.
\end{proof}

\section{Proof of Corollaries}
\label{s_corollaries}

\begin{proof}[Proof of \eqref{c_main}] Since $X$ is covered by rational curves which 
intersect an irreducible component of $V$, the main theorem of \cite{Campana91} implies
that the image of the fundamental group of $V$ is of finite index in the fundamental group
of $X$.

 Let $\pi\colon\map Y.X.$ be any \'etale morphism of degree $r$, where $Y$ is connected,
let $W$ be the inverse image of $V$, and set $\Gamma=\pi^*\Delta$.  Then
$K_Y+\Gamma=\pi^*(K_X+\Delta)$, so that $W$ is the locus of log canonical singularities of
$K_Y+\Gamma$ and $-(K_Y+\Gamma)$ is big and nef.  Note that $W$ is connected, by the
connectedness theorem.

The proof now divides into two cases.  If $V$ is non-empty, then the result follows
immediately from the fact that $\pi$ is arbitrary and $W$ is always connected.  

Now suppose that $V$ is empty.  In this case, we adopt the convention that the fundamental
group of the empty set is the trivial group.  By what we have already proved the
fundamental group is finite and so it suffices to prove that the algebraic fundamental
group is trivial.  Now
$$
\chi(\ring Y.)=r\chi(\ring X.).
$$
On the other hand, 
$$
h^i(Y,\ring Y.)=h^i(X,\ring X.)=0 \qquad \text{for $i>0$,} \qquad 
$$
by Kawamata-Viehweg vanishing, so that 
$$
\chi(\ring Y.)=\chi(\ring X.)=1,
$$
whence $r=1$.  But then the fundamental group of $X$ is trivial and the fundamental
group of $V$ certainly surjects onto the fundamental group of $X$.  
\end{proof}

\begin{lemma}\label{l_semi} Let $(X,\Delta)$ be a kawamata log terminal pair and
let $f\colon\map X.S.$ be a projective morphism.  Suppose that $-(K_X+\Delta)$ is
relatively nef and $-K_X$ is relatively big.

 Then $-(K_X+\Delta)$ is relatively semiample. 
\end{lemma}
\begin{proof} It suffices to exhibit a divisor $\Theta$, such that $K_X+\Theta$ is 
kawamata log terminal and $-(K_X+\Theta)$ is ample, since then, by the
base point free theorem, every nef divisor is relatively semiample.

  By assumption we may write 
$$
-K_X\sim_{\mathbb{Q},f} A+B,
$$
where $A$ is relatively ample and $B$ is effective.  Let 
$$
\Theta=(1-\epsilon)\Delta+\epsilon B.
$$

 Then 
$$
-(K_X+\Theta)\sim_{\mathbb{Q},f} -(1-\epsilon)(K_X+\Delta)+\epsilon A
$$
is relatively ample, as it is the sum of a relatively nef divisor and a relatively
ample divisor.  In particular $K_X+\Theta$ is $\mathbb{Q}$-Cartier.  But then
$$
K_X+\Theta=K_X+(1-\epsilon)\Delta+\epsilon B,
$$
is certainly kawamata log terminal, for $\epsilon$ small enough.  
\end{proof}

\begin{proof}[Proof of \eqref{c_rcck}] By \eqref{l_semi} $\ring X.(-m(K_X+\Delta))$ is 
relatively generated, for some $m>0$, and so the result follows from \eqref{t_main}, since
the locus of log canonical singularities is empty.
\end{proof}

\begin{proof}[Proof of \eqref{c_rccl}] By \eqref{t_main}, we already know that the fibres
of $f$ are rationally chain connected modulo the locus $V$ of log canonical singularities.
It remains to prove that the fibres of $V$ over $S$ are rationally chain connected.  By
induction on the dimension, it suffices to prove that if $W$ is a log canonical centre,
then $W$ is rationally chain connected modulo the union $R$ of the log canonical centres
properly contained in $W$.  As in the proof of \eqref{t_main}, by \eqref{l_alter}, we may
assume that $\Delta$ contains an ample divisor and that some multiple of $K_X+\Delta$
is linearly equivalent to zero.  Let $g\colon\map Y.X.$ be a log resolution of the pair
$(X,\Delta)$.  We may write
$$
K_Y+\Gamma=g^*(K_X+\Delta)+E,
$$
where $\Gamma$ and $E$ are effective, with no common components and $E$ is
$g$-exceptional.  Let $F$ a log canonical centre of $K_Y+\Gamma$ which dominates $W$,
minimal with this property.  Then $F$ is the intersection of divisors of log discrepancy
zero, with centre $W$.  In particular we may find a divisor $\Theta$ such that
$$
(K_Y+\Gamma)|_F=K_F+\Theta.
$$
We are going to apply \eqref{c_uni} to $h\colon\map F.W.$.  Let $t\colon\rmap W.Z.$ be
a rational map.  Blowing up further, we may assume that the induced map $\psi\colon\map
F.Z.$ is a morphism.  Note that the image in $W$ of the locus of log canonical
singularities of $K_F+\Theta$ is equal to $R$, by minimality of $F$.  So we may assume
that that the locus of log canonical singularities of $K_F+\Theta$ does not dominate $Z$.
We need to check that conditions (2-4) hold.  (2) holds, exactly as in the proof of
\eqref{t_main}.  (3) holds by (4) of \eqref{t_fibres} and (4) holds by assumption.
\end{proof}

\begin{proof}[Proof of \eqref{c_dlt}] By (2.43) of \cite{KM98}, we may assume that 
the pair $(X,\Delta)$ is kawamata log terminal, and in this case we may apply
\eqref{t_main} with $f$ the identity.  \end{proof}

\begin{proof}[Proof of \eqref{c_covered}] Immediate from \eqref{c_dlt}.  
\end{proof}

\begin{proof}[Proof of \eqref{c_equivalent}] If $X$ is rationally connected then it is
certainly rationally chain connected, so suppose that $X$ is rationally chain connected.
Let $g\colon\map Y.X.$ be a resolution of singularities.  By \eqref{c_dlt} $g$ has
rationally chain connected fibres.  It follows that $Y$ is rationally chain connected. 
But then $Y$ is rationally connected, as it is smooth, and so $X$ is rationally connected 
as well.  
\end{proof}

\begin{proof}[Proof of \eqref{c_over}] We work locally in a neighbourhood of the generic
point of $T$.  \eqref{l_semi} implies that $-(K_X+\Delta)$ is relatively semiample.  As in
the proof of \eqref{t_main} and \eqref{l_alter}, we may therefore assume $-(K_X+\Delta)$
is ample and that $\Delta$ contains an ample divisor $A$,

Suppose that $T\neq S$.  Pick any ample divisor $H$ in $S$ which contains $T$.  Let $t$ be
the largest rational number such that $K_X+\Delta+tf^*H$ is log canonical.  Possibly
perturbing $\Delta$ by an ample divisor in $X$, we may assume that there is a unique log
canonical centre $V$ of $K_X+\Delta+tf^*H$.  By connectedness $V$ has connected fibres
over $T$.  Since $V$ is the unique log canonical centre, it is certainly a minimal log
canonical centre.  Thus $V$ is normal.  It follows that $V$ has irreducible fibres.  By
Kawamata's subadjunction formula, see Theorem 1 of \cite{Kawamata97}, there is a boundary
$\Theta$ on $V$, such that
$$
(K_X+\Delta)|_V=K_V+\Theta,
$$
where $K_V+\Theta$ is kawamata log terminal.  Possibly replacing $A$ by a linearly 
equivalent divisor, we may assume that $\Theta$ contains an ample divisor.  Thus we 
are free to restrict to $V$ and replace $X$ by $V$ and $S$ by the image of $V$.  

Thus by induction on the dimension, we may assume that $T=S$.  In this case the general
fibre over $T$ is kawamata log terminal and so by \eqref{t_main} and \eqref{c_equivalent}
the general fibre is rationally connected.  \end{proof}

\begin{proof}[Proof of \eqref{c_section}] Pick a resolution $g\colon\map Y.X.$ and let 
$\pi\colon\map Y.S.$ be the composition.  Pick a curve $C\subset S$, and let 
$E$ be a subvariety of $Y$ with rationally connected fibres over $C$, whose existence
is guaranteed by \eqref{c_over}.   Now apply the main result of \cite{GHS03}.  \end{proof}

\begin{proof}[Proof of \eqref{c_chow}] Even though this would seem to be a fairly standard
result, we could not find a reference, and so we include a proof for completeness.

First note that if $\pi$ is any surjective morphism, then $\pi_*$ is always surjective, as
$\ch^0(S)$ is generated by the points of $S$.  Now suppose that $\alpha$ is a cycle in
$\ch^0(Y)$, whose image $\beta$ in $\ch^0(S)$ is zero.  It remains to prove that $\alpha$
is equivalent to zero.  First pick a curve $\Sigma$ in $Y$ which is the intersection of
general hyperplanes, containing the support of $\alpha$.  Then we may replace $\alpha$ by
an equivalent cycle in $\Sigma$, in such a way that the support of $\alpha$ belongs to the
general fibre.

As $\beta$ is equivalent to zero, by definition of $\ch^0(S)$, it follows that there is a
curve $C$ in $S$ such that the support of $\beta$ belongs to $C$, and such that $\beta$ is
linearly equivalent to zero in $C$.  By \eqref{c_section} we may find a curve $D\subset Y$
which maps birationally down to $C$.  Let $\alpha'$ be the pullback of $\beta$ to $D$.  As
$\alpha$ is equivalent to zero in $C$, $\alpha'$ is certainly equivalent to zero in $D$.  

On the the other hand, by our choice of $\alpha$, the fibres of $\pi$ over the support of
$\beta$ are rationally connected.  It follows that $\alpha$ and $\alpha'$ are equivalent,
so that $\alpha$ is indeed zero in $\ch^0(Y)$.  \end{proof}

\begin{proof}[Proof of \eqref{c_iitaka}] If $K_X+\Delta$ is kawamata log terminal, then 
the locus of log canonical singularities is empty and the result is therefore immediate
from \eqref{c_uni}.  \end{proof}

\begin{proof}[Proof of \eqref{c_fano}] By assumption the Iitaka fibration of $-(K_X+\Delta)$ 
is birational, and as rational connectedness is a birational invariant, the result follows
from \eqref{c_iitaka}.
\end{proof}

\begin{proof}[Proof of \eqref{c_simply}] The pair $(X,\Delta)$ is kawamata log terminal iff 
the locus of log canonical singularities is empty, and so this result is an immediate
consequence of \eqref{c_main}.  \end{proof}

\bibliographystyle{hamsplain} 

\bibliography{/home/mckernan/Jewel/Tex/math}

\end{document}